\documentclass[american,parskip=half,12pt]{scrartcl}

\usepackage[T1]{fontenc}
\usepackage[utf8]{inputenc}
\usepackage{babel}
\usepackage{csquotes}
\usepackage[shortlabels]{enumitem}

\usepackage{amssymb,amsmath,bm,mathtools}
\usepackage{mathabx} % odiv
\usepackage[hyperref,amsmath,thmmarks]{ntheorem}

\usepackage{graphicx}
\usepackage[justification=raggedright]{subcaption}
\usepackage{xcolor}
\usepackage{pgfplots}
\graphicspath{{../images/},{../matlab/images/}}
\usepgfplotslibrary{external}

\usepackage[english,pdfusetitle,colorlinks=true,linkcolor=blue,citecolor=green!50!black,urlcolor=blue,bookmarks=true]{hyperref}

\usepackage[ruled]{algorithm2e}

\usepackage{todonotes}
\usepackage{aliascnt}
\usepackage{siunitx}

\newaliascnt{proposition}{lemma}

\aliascntresetthe{proposition}

\newaliascnt{corollary}{lemma}

\aliascntresetthe{corollary}

\newaliascnt{theorem}{lemma}
\newtheorem{theorem}[theorem]{Theorem}
\aliascntresetthe{theorem}

\theorembodyfont{\normalfont}
\newaliascnt{definition}{lemma}

\aliascntresetthe{definition}

\newaliascnt{assumption}{lemma}

\aliascntresetthe{assumption}

\newaliascnt{notation}{lemma}

\aliascntresetthe{notation}

\newaliascnt{example}{lemma}

\aliascntresetthe{example}

\newaliascnt{experiment}{lemma}

\aliascntresetthe{experiment}

\newaliascnt{remark}{lemma}

\aliascntresetthe{remark}

\theoremstyle{nonumberplain}
\theoremseparator{:}
\theoremheaderfont{\normalfont\itshape}

\theoremsymbol{\ensuremath{\square}}

\newcommand{\Z}{\ensuremath{\mathbb{Z}}}

\newcommand{\R}{\ensuremath{\mathbb{R}}}
\newcommand{\C}{\ensuremath{\mathbb{C}}}

\newcommand{\SO}{\ensuremath{\mathrm{SO}(3)}}

\newcommand{\abs}[1]{\ensuremath{\left\vert#1\right\vert}}

\newcommand{\e}{\mathrm{e}}

\renewcommand{\i}{\mathrm{i}}

\DeclareMathOperator*{\argmin}{arg\,min}

\DeclareMathOperator*{\sgn}{sgn}

\DeclareMathOperator{\grad}{grad}

\newcommand{\dd}{\, \mathrm{d} }
\newcommand{\bn}{{\bm{n}}}

\newcommand{\bl}{{\bm \ell}}

\newcommand{\bx}{{\bm x}}

\newcommand{\bk}{{\bm k}}

\newcommand{\ui}{u^{\mathrm{inc}}}

\newcommand{\utot}{u^{\mathrm{tot}}}

\newcommand{\rs}{r_{\mathrm{s}}}

\newcommand{\rM}{r_{\mathrm{M}}}

\title{Total Variation-Based Reconstruction and Phase Retrieval for Diffraction Tomography with an Arbitrarily Moving Object}

\date{}
\author{%
  Robert Beinert\thanks{TU Berlin, Institute of
    Mathematics, MA 4-3, Straße des 17. Juni 136, 10623 Berlin, Germany.
  }\\%
  {\footnotesize\href{mailto:beinert@math.tu-berlin.de}{beinert@math.tu-berlin.de}}
  \and 
  Michael Quellmalz\footnotemark[1]\\%
  {\footnotesize\href{mailto:quellmalz@math.tu-berlin.de}{quellmalz@math.tu-berlin.de}}
}

\begin{document}

\def\sectionautorefname{Section}
\def\subsectionautorefname{Section} % works only after begin document

\maketitle

\begin{abstract}
  \noindent
  We consider the imaging problem of the reconstruction of a three-di\-men\-sio\-nal object 
  via optical diffraction tomography under the assumptions of the Born approximation.
  Our focus lies in the situation that a rigid object performs an irregular, time-dependent rotation under acoustical or optical forces.
  In this study, we compare reconstruction algorithm in case
  i) that two-dimensional images of the complex-valued wave are known, or
  ii) that only the intensity (absolute value) of these images can be measured, which is the case in many practical setups.
  The latter phase-retrieval problem can be solved by an all-at-once approach based 
  utilizing a hybrid input-output scheme with TV regularization.
\end{abstract}

\section{Introduction}

We are interested in the tomographic reconstruction of a three-dimensional (3D) rigid object, which is illuminated from various directions. 
In \emph{optical diffraction tomography}, 
the wavelength of the imaging wave, visible light with wavelength hundreds of nanometers, is in a similar order of magnitude as features of the \textmu{}m-sized object we want to reconstruct.
This setup substantially differs from the X-ray computerized tomography,
where the wavelength is much shorter and therefore the light can be assumed to propagate along straight lines.

\begin{figure}
  \centering
  \begin{tikzpicture}[scale=0.65]
    % axes
    \draw[->] (-5,0)--(5,0) node[right] {$x_1$};
    \draw[->] (0,-4.5)--(0,4.5) node[right] {$x_3$};
    \draw[->,opacity=.7] (-2,-2)--(2,2) node[right] {$x_2$};
    % object
    \fill[red,opacity=0.3] (0,0) circle (2);
    % \node at (1,0.5) {\footnotesize{object}};
    \node at (1,-0.8) {$f$};
    \draw[->]  (-2.46, 0.43) arc (170:100:2.5);
    % incident wave
    \draw[dashed] (-5,-3.5) -- (5,-3.5);
    \draw[dashed] (-5,-4) -- (5,-4);
    \draw[dashed] (-5,-4.5) -- (5,-4.5);
    \node at (-7,-4) {\footnotesize{incident field}};
    \draw[->] (5.5,-4.5) -- (5.5,-3.5) node[right] {$u^{\mathrm{inc}}$};
    % measurement plane
    \draw[line width = 3pt] (5,3.5) -- (-5,3.5) node [left,align=center] {\footnotesize{measurement plane}};
    \node at (7,3.5) {$x_3 = r_{\mathrm M}$};
  \end{tikzpicture} 
  \caption{Experimental setup with measurement plane located at $x_3=r_{\mathrm M}$, see \cite{BeiQue22}.
  }	\label{fig:trans}
\end{figure}
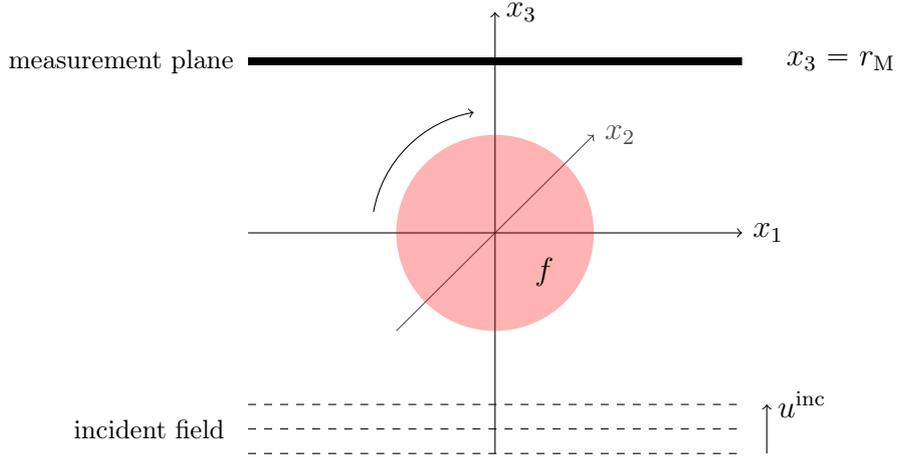

Formally, in optical diffraction tomography we want to recover an object or, more specifically, its \emph{scattering potential}
\begin{equation*}
  f\colon\mathbb R^3\to\mathbb R,
\end{equation*}
which has compact support.
At each time step $t\in[0,T]$, 
the object undergoes a rotation described by a \emph{rotation matrix} $R_t\in SO(3)$,
such that its scattering potential becomes $f(R_t\cdot)$.
Here, the rotation is not restricted to some direction, we only require that $R_t$ depends continuously differentiably on $t$,
which is the case for objects in acoustical traps, cf.\ \cite{ElbSchm20}.
We further assume that the rotation $R_t$ is known beforehand;
for the detection of motion from the tomographic data we refer to \cite{QueElbSchSte22}.
The object is illuminated by an \emph{incident, plane wave}
\begin{equation*}
  u^{\mathrm{inc}}(\boldsymbol x) := \mathrm e^{\mathrm i k_0 x_3}
  ,\quad \boldsymbol x = (x_1,x_2,x_3)\in\R^3,
\end{equation*}
of \emph{wave number} $k_0>0$,
which induces a \emph{scattered wave} $u_t\colon\mathbb R^3\to\mathbb C$.
We measure the resulting \emph{total field} 
\begin{equation*}
  \utot_t(\bx) = u_t(\bx)+u^{\mathrm{inc}}(\bx),
  \quad \bx \in \R^3,
\end{equation*}
at a fixed measurement plane $x_3 = r_{\mathrm M}$,
see Fig.~\ref{fig:trans}.

Under the \emph{Born approximation},
the scattered wave $u_t$ is a solution of the partial
differential equation
\begin{equation*}
  -(\Delta + k_0^2)u_t = (u_t+u^{\mathrm{inc}})\, f(R_t \cdot),
\end{equation*}
see \cite{NatWue01}.
We assume the \emph{Born approximation} to be valid,
which is the case for small, mildly scattering objects, see \cite[§\,3.3]{FauKirQueSchSet21}.
Then the relationship between the measured field $u_t$ and the desired function $f$ can be expressed via the \emph{Fourier diffraction theorem} \cite{KirQueRitSchSet21,KakSla} by
\begin{equation} \label{eq:recon}
  \sqrt{\frac{2}{\pi}} \,
  \frac{\kappa}{\mathrm i \mathrm e^{\mathrm i\kappa \rM}} \, 
  \mathcal F_{1,2}[u_t](k_1,k_2, \rM)
  =
  \mathcal F[f]\left(R_t (k_1,k_2,\kappa-k_0) \right),
  \quad k_1^2+k_2^2 < k_0^2,
\end{equation}
where
\begin{equation*}
  \kappa(k_1,k_2) := \sqrt{k_0^2-k_1^2-k_2^2},
  \quad k_1^2+k_2^2 < k_0^2.
\end{equation*}
Here the 3D \emph{Fourier transform} $\mathcal F$
and the \emph{partial Fourier transform} $\mathcal F_{1,2}$ in the first two coordinates
are given by
\begin{equation*}
  \mathcal F [f](\bk)
  \coloneqq (2\pi)^{-\frac32} \int_{\R^3} f(\bx)\,\e^{-\i\, \bx\cdot\bk} \dd \bx
\end{equation*}
and
\begin{equation*}
  \mathcal F_{1,2} [f](\bk',x_3)
  \coloneqq (2\pi)^{-1} \int_{\R^{2}} f(\bx',x_3)\, \e^{-\i\, \bx'\cdot\bk'} \dd \bx'
  \vspace{5pt}
\end{equation*}
for $\bk \in \R^3$ and $\bk' \in \R^2$.
The left-hand side of \eqref{eq:recon} is fully determined by the measurements $\utot_t(\cdot,\cdot,\rM)$,
and the right-hand side provides non-uniform samples of the Fourier transform $\mathcal F[f]$,
evaluated on the union of semispheres which contain the origin and are rotated by $R_t$.
Therefore, the reconstruction problem in diffraction tomography can be seen as a problem of inverting the 3D Fourier transform $\mathcal F$ given non-uniformly sampled data in the so-called \emph{$k$-space}.
If these samples have a positive Lebesgue measure,
the scattering potential is uniquely defined by the measurements.

\begin{theorem}[Unique inversion, {\cite[Thm.~3.2]{BeiQue22}}] \label{thm:unique}
Let $f\in L^1(\R^3)$ have a compact support in $\{\bx \in \R^3 : \|\bx\|_2 \le \rM \}$,
and let the  Lebesgue measure of
\begin{equation}
  \label{eq:Y}
  \mathcal Y 
  \coloneqq 
  \{ R_t (k_1,k_2,\kappa-k_0) : k_1^2 + k_2^2 \le k_0^2,\, t\in [0,T]\}
  \subset \R^3
\end{equation}
be positive.  Then $f$ is uniquely determined by
$\utot_t(\cdot,\cdot,\rM)$, $t \in [0,T]$.
\end{theorem}

In §\,\ref{sec:rec} of this paper, we consider algorithms for the tomographic reconstruction of $f$ from $\utot_t(\cdot,\cdot,\rM)$, $t\in[0,T]$.
Since in practice one can often measure only intensities $|\utot_t|$, cf.\ \cite{MalDev93},
in §\,\ref{sec:pr} we will incorporate phase retrieval methods to retrieve $f$.
In contrast to \cite{BeiQue22}, we focus here on the numerical studies for the case of an arbitrarily moving rotation axis.

\section{Reconstruction methods}
\label{sec:rec}

We consider three different reconstruction methods for the situation that we know the complex-valued field $u_t(\cdot,\cdot,\rM)$.
First, the \emph{filtered backpropagation} is based on a truncation and discretization of the inverse 3D Fourier transform,
\[
f_{\mathrm{bp}}(\boldsymbol x)
=
\int_{\mathcal Y}
\mathcal F[f](\boldsymbol y)
\,\mathrm e^{\mathrm i \boldsymbol x\cdot\boldsymbol y}
\,\mathrm d\boldsymbol y,
\quad \bx\in\R^3,
\]
where $\mathcal Y$ is the set where the data $\mathcal Ff$ is available in $k$-space.
We state the following filtered backpropagation formula for a moving rotation axis, the special case for constant axis is well established \cite{Dev82}.

\begin{theorem}[Filtered Backpropagation, {\cite[Thm.~4.1]{KirQueRitSchSet21}}]
  \label{thm:reconstruction}
  Let the assumptions of the Fourier diffraction theorem be satisfied.
  Let $\alpha \in C^1([0,T],\R)$ and $\boldsymbol n\in C^1([0,T],\mathbb S^2)$. 
  For $t\in[0,T]$, we denote the rotation around the axis $\boldsymbol n(t)$ with the angle $\alpha(t)$ by
  \begin{equation*}
    R_{\boldsymbol n(t),\alpha(t)}^\top \coloneqq 
    \begin{pmatrix}
      n_1^2(1-\mathrm{c}) + \mathrm{c} & n_1 n_2 (1-\mathrm{c}) - n_3 \mathrm{s} & 
      n_1n_3(1-\mathrm{c}) + n_2 \mathrm{s}\\
      n_1n_2(1-\mathrm{c}) + n_3 \mathrm{s} & n_2^2 (1-\mathrm{c}) + \mathrm{c} & 
      n_2n_3(1-\mathrm{c}) - n_1 \mathrm{s}\\
      n_1n_3(1-\mathrm{c}) - n_2\mathrm{s} & n_2 n_3 (1-\mathrm{c}) + n_1 \mathrm{s} & 
      n_3^2(1-\mathrm{c}) + \mathrm{c}
    \end{pmatrix}
    \in \SO,
  \end{equation*}
  where $\boldsymbol n(t)=(n_1,n_2,n_3)^\top$, $\mathrm{c} \coloneqq \cos \alpha(t)$ and $\mathrm{s} \coloneqq \sin \alpha(t) $.
  Denote by $u_t$ the scattered wave according to the rotated object $f \circ R_{\boldsymbol n(t),\alpha(t)}$.
  Next, let 
  \begin{equation*}
    \mathcal{U} \coloneqq \{(k_1,k_2,t) \in \mathbb R^3 : k_1^2 + k_2^2 < k_0^2, \, 0\le t \le T \}
  \end{equation*}
  and
  \begin{equation*}\label{eq:T}
    \Phi: \mathcal{U} \to \mathbb{R}^3, \quad T(k_1,k_2,t) \coloneqq R_{\boldsymbol n(t),\alpha(t)}(k_1,k_2,  \kappa - k_0)^\top
  \end{equation*}
  be the map that covers the accessible domain $\mathcal Y = \Phi(\mathcal U)$ in k-space.
  Then, for all $\bx \in \mathbb R^3,$
  \begin{equation}\label{eq:reconstruction}
    f_{\mathrm{bp}}(\bx) = \frac{-\mathrm i}{2 \pi^2} \int_{\mathcal{U}} \mathrm e^{\mathrm i\, \Phi(k_1,k_2,t) \cdot \bx}\, \mathcal F_{1,2}[u_t](k_1,k_2, \rM) \,\frac{\kappa \mathrm e^{-\mathrm i \kappa \rM} |{\nabla \Phi(k_1,k_2,t)}|}{\operatorname{Card}(\Phi^{-1}(\Phi(k_1,k_2,t)))} \,\mathrm d(k_1,k_2,t),
  \end{equation}
  where $\operatorname{Card}$ denotes the cardinality of a set 
  and $|\nabla \Phi|$ is the absolute value of the Jacobian determinant of $\Phi$.
\end{theorem}

If $\bn(t)$ and $\alpha(t)$ are known analytically in \eqref{eq:reconstruction}, then one can compute $|\nabla \Phi|$,
but the Banach indicatrix or Crofton symbol $\operatorname{Card}(\Phi^{-1}(\Phi(k_1,k_2,t)))$ is hard to determine analytically.
However, since the backpropagation tends to produce artifacts due to the limited coverage $\mathcal Y$ in $k$-space, we consider alternative approaches.
The second approach consists in finding 
\begin{equation} \label{eq:ndft}
\argmin_{\boldsymbol f}
\| \boldsymbol F \boldsymbol f - \boldsymbol g \|^2,
\end{equation}
where $\boldsymbol F$ denotes a discretization of $\mathcal F$,
a \emph{nonuniform discrete Fourier transform} (NDFT) \cite{PlPoStTa18},
$\boldsymbol f$ are uniform samples of $f$ and $\boldsymbol g$ is the uniformly sampled data on the left-hand side of \eqref{eq:recon}, which is determined by the measurements.
The resulting normal equation can be solved with a \emph{conjugate gradient} (CG) method,
making use of fast implementations \cite{KeKuPo09} of the NDFT and its adjoint.
Such methods have also been applied in X-ray imaging \cite{PoSt01IMA}, magnetic resonance
imaging \cite{KnKuPo}, and spherical tomography \cite{HiQu15,HiQu16}.

Since the CG reconstruction is sensitive regarding to larger measurement errors
in order to be applied in phase retrieval,
we further add a \emph{total variation} (TV) regularization term
and incorporate the non-negativeness of the object, see \cite[§\,5.3]{BeiQue22}.
For the discretized object
\begin{equation*}
  (\bm f)_{\bl} \coloneqq f(\bx_\bl),
  \quad
  \bx_\bl \coloneqq \frac{1}{L} \, \bl \in \R^3,
  \quad
  \bl \in \mathcal I,
\end{equation*}
on a uniform Cartesian grid,
where $\mathcal I \subset \Z^3$ denotes the employed index set,
the \emph{total variation} is given by
\begin{equation*}
  \mathrm{TV}(\bm f)
  \coloneqq \sum_{\bl \in \mathcal I}
  \| ( \grad \bm f )_{\bl} \|_2,
\end{equation*}
see e.g.\ \cite{BreLor18},
where the discrete gradient $\grad$ may be computed using finite forward differences.
To incorporate the non-negativity constraint,
we exploit the characteristic function
\begin{equation*}
  \chi_{{\ge 0}}(\bm f)
  \coloneqq
  \begin{cases}
    0, & \text{if} \; f_{\bl} \ge 0 \; \text{for all} \; \bl \in
    \mathcal I, \\
    +\infty, & \text{otherwise}.
  \end{cases}
\end{equation*}
Instead of minimizing \eqref{eq:ndft},
we propose to invert the NDFT by computing
\begin{equation}
  \label{eq:var-prob}
  \argmin_{\bm f} \quad
  \chi_{{\ge0}}(\bm f)
  + \frac12 \, \|\bm F
    (\bm f) -
    \bm g\|^2_{2} +
  \lambda \, \mathrm{TV}(\bm f).
\end{equation}
The minimization problem \eqref{eq:var-prob} can be solved
using a \emph{primal-dual method} \cite{ChaCasCreNovPoc10} with backtracking \cite{YokHon17}.
We compare the different reconstruction methods numerically for the rotation axis 
\begin{equation}
  \label{eq:axis}
  \bn(t)
  =
  \left(\sin\bigl(\tfrac12\pi \sin t\bigr),
    \cos\bigl(\tfrac12\pi \sin t\bigr),
    0\right)^\top
\end{equation}
and angle
$\alpha(t) = t$,
$t\in[0,2\pi]$,
for which the backpropagation formula is given in \cite[Ex.~4.8]{KirQueRitSchSet21}.
Our forward data $\utot_t$ is created based on the Born approximation.
In Fig.~\ref{fig:cell3d-rec}, we see that the filtered backpropagation falls behind the other methods.
Adding 5\,\% Gaussian noise to the data $u_t+\ui$, we see in Fig.~\ref{fig:cell3d-rec-noise} that the TV-regularized inverse is superior to the others.
Here, we perform 20 iterations for the CG and 200 iterations for the primal-dual TV-regularization.
The forward model to generate the data also uses Born approximation, but a different algorithm based on convolution.

\tikzset{font=\tiny}
\newcommand{\modelwidth} {3.5cm}
\pgfplotsset{
  colormap={parula}{
    rgb255=(53,42,135)
    rgb255=(15,92,221)
    rgb255=(18,125,216)
    rgb255=(7,156,207)
    rgb255=(21,177,180)
    rgb255=(89,189,140)
    rgb255=(165,190,107)
    rgb255=(225,185,82)
    rgb255=(252,206,46)
    rgb255=(249,251,14)}}
% https://tex.stackexchange.com/questions/247940/parula-colormap-in-pgfplots
\pgfmathsetmacro{\xmin}{-42.426407}
\pgfmathsetmacro{\xmax}{42.072853}

\begin{figure}[p]\centering
  \pgfmathsetmacro{\cmin} {-0.05} \pgfmathsetmacro{\cmax} {1.05}
  \begin{tikzpicture}
    \begin{axis}[
      width=\modelwidth,
      height=\modelwidth,
      enlargelimits=false,
      scale only axis,
      axis on top,
      colorbar,colorbar style={
        width=.15cm, xshift=-0.9em, 
        point meta min=\cmin,point meta max=\cmax,
      },
      ]
      \addplot graphics [
      xmin=\xmin, xmax=\xmax,  ymin=\xmin, ymax=\xmax,
      ] {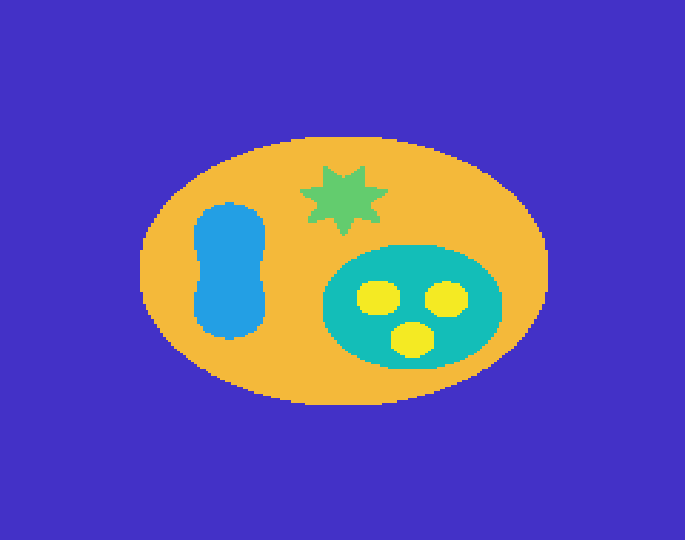};
    \end{axis}
  \end{tikzpicture}
  \qquad
  \begin{tikzpicture}
    \begin{axis}[
      width=\modelwidth,
      height=\modelwidth,
      enlargelimits=false,
      scale only axis,
      axis on top,
      colorbar,colorbar style={
        width=.15cm, xshift=-0.9em, 
        point meta min=\cmin,point meta max=\cmax,
      },
      ]
      \addplot graphics [
      xmin=\xmin, xmax=\xmax,  ymin=\xmin, ymax=\xmax,
      ] {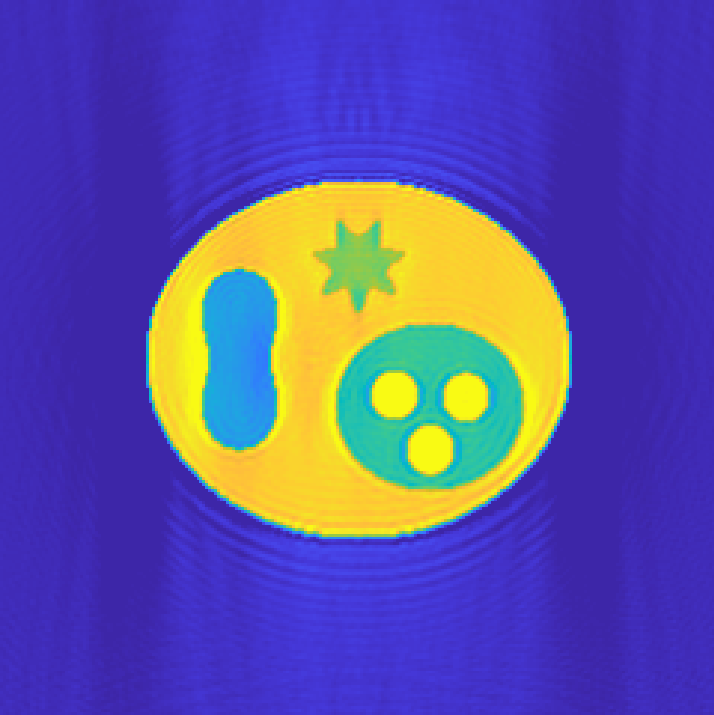};
    \end{axis}
  \end{tikzpicture}\\[10pt]
  \begin{tikzpicture}
    \begin{axis}[
      width=\modelwidth,
      height=\modelwidth,
      enlargelimits=false,
      scale only axis,
      axis on top,
      colorbar,colorbar style={
        width=.15cm, xshift=-0.9em, 
        point meta min=\cmin,point meta max=\cmax,
      },
      ]
      \addplot graphics [
      xmin=\xmin, xmax=\xmax,  ymin=\xmin, ymax=\xmax,
      ] {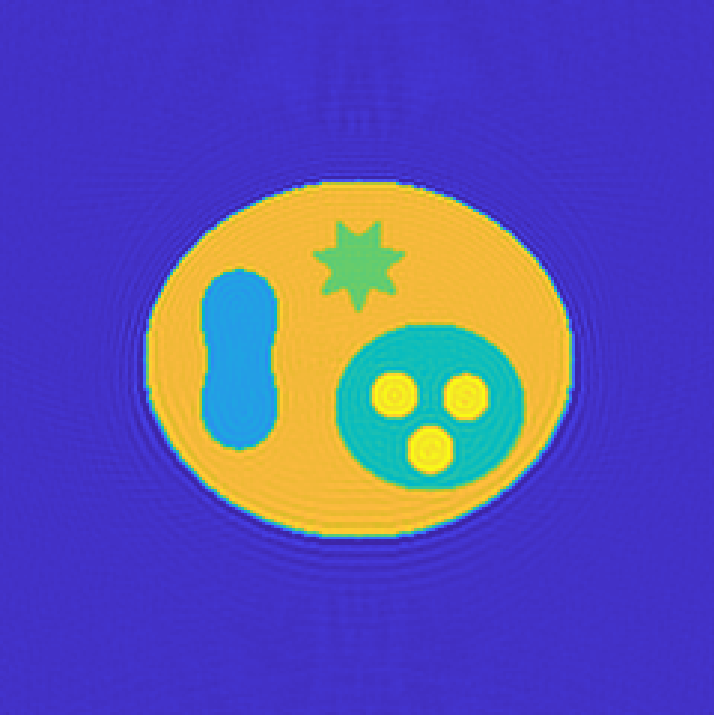};
    \end{axis}
  \end{tikzpicture}
  \qquad
  \begin{tikzpicture}
    \begin{axis}[
      width=\modelwidth,
      height=\modelwidth,
      enlargelimits=false,
      scale only axis,
      axis on top,
      colorbar,colorbar style={
        width=.15cm, xshift=-0.9em, 
        point meta min=\cmin,point meta max=\cmax,
      },
      ]
      \addplot graphics [
      xmin=\xmin, xmax=\xmax,  ymin=\xmin, ymax=\xmax,
      ] {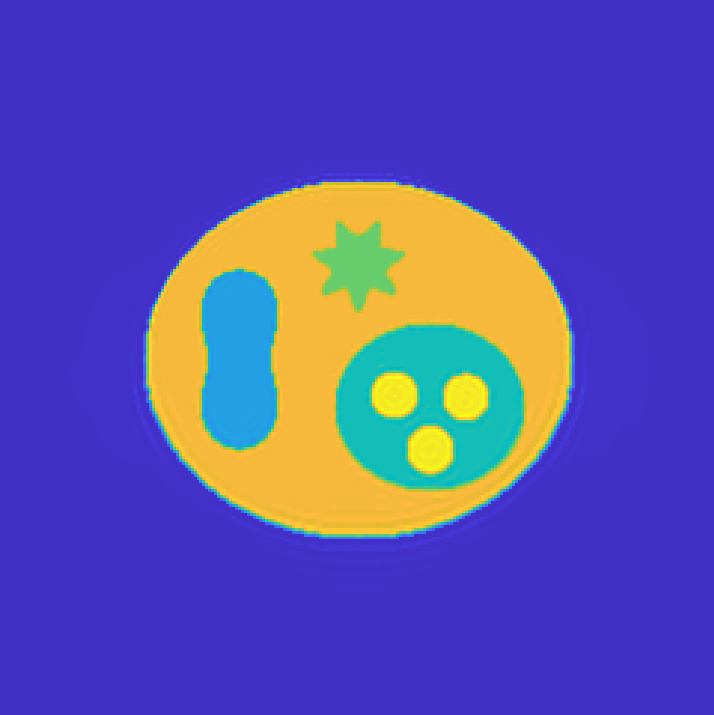};
    \end{axis}
  \end{tikzpicture}
  \caption{Reconstructions with moving axis (exact data).
    Top left: Ground truth of 3D phantom $f$.
    Top right: Filtered backpropagation (PSNR 28.55, SSIM 0.781).
    Bottom left: CG Reconstruction (PSNR 33.75, SSIM 0.962). 
    Bottom right: Primal-dual with TV and regularization parameter $\lambda=0.01$ (PSNR 34.54, SSIM 0.991). 
    \label{fig:cell3d-rec}}
\end{figure}

\begin{figure}[p]\centering
\pgfmathsetmacro{\cmin} {-0.05} \pgfmathsetmacro{\cmax} {1.05}
\begin{tikzpicture}
  \begin{axis}[
    width=\modelwidth,
    height=\modelwidth,
    enlargelimits=false,
    scale only axis,
    axis on top,
    colorbar,colorbar style={
      width=.15cm, xshift=-0.9em, 
      point meta min=\cmin,point meta max=\cmax,
    },
    ]
    \addplot graphics [
    xmin=\xmin, xmax=\xmax,  ymin=\xmin, ymax=\xmax,
    ] {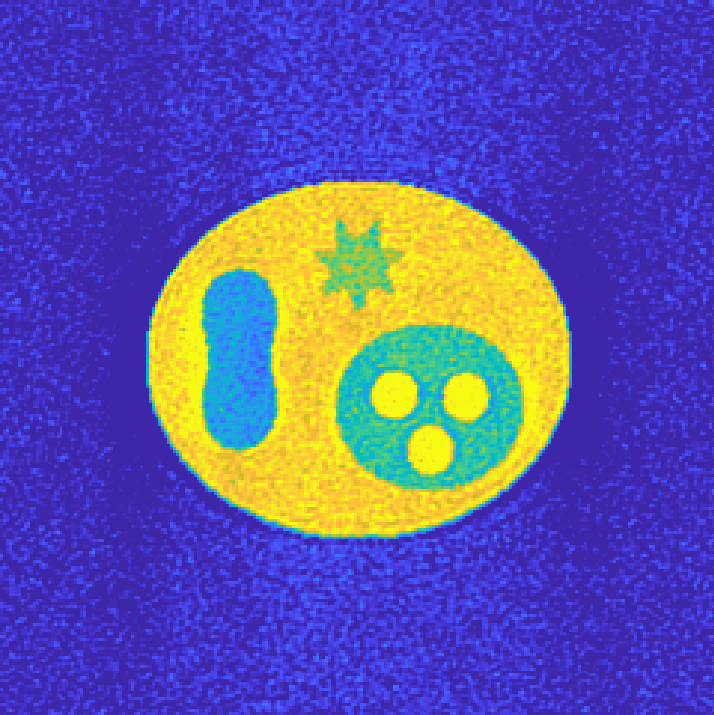};
  \end{axis}
\end{tikzpicture}
\begin{tikzpicture}
  \begin{axis}[
    width=\modelwidth,
    height=\modelwidth,
    enlargelimits=false,
    scale only axis,
    axis on top,
    colorbar,colorbar style={
      width=.15cm, xshift=-0.9em, 
      point meta min=\cmin,point meta max=\cmax,
    },
    ]
    \addplot graphics [
    xmin=\xmin, xmax=\xmax,  ymin=\xmin, ymax=\xmax,
    ] {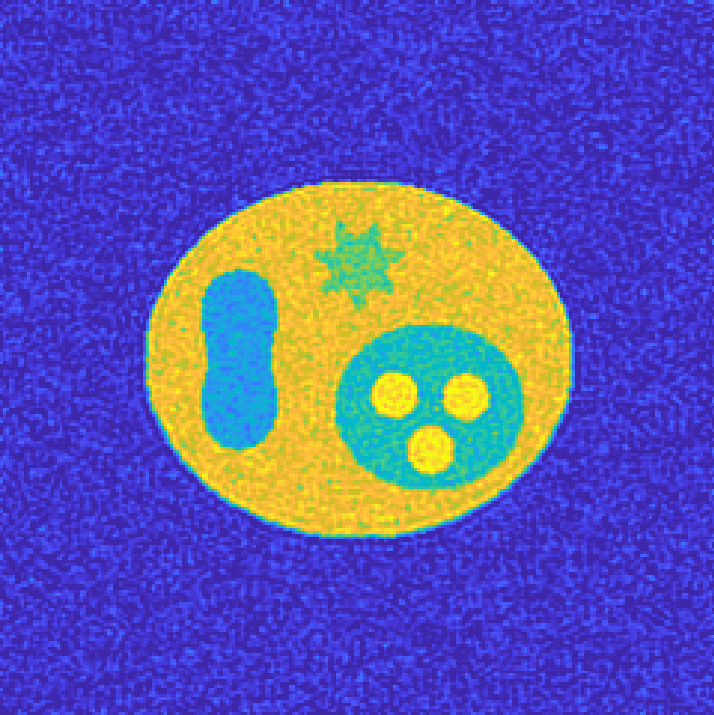};
  \end{axis}
\end{tikzpicture}
\begin{tikzpicture}
  \begin{axis}[
    width=\modelwidth,
    height=\modelwidth,
    enlargelimits=false,
    scale only axis,
    axis on top,
    colorbar,colorbar style={
      width=.15cm, xshift=-0.9em, 
      point meta min=\cmin,point meta max=\cmax,
    },
    ]
    \addplot graphics [
    xmin=\xmin, xmax=\xmax,  ymin=\xmin, ymax=\xmax,
    ] {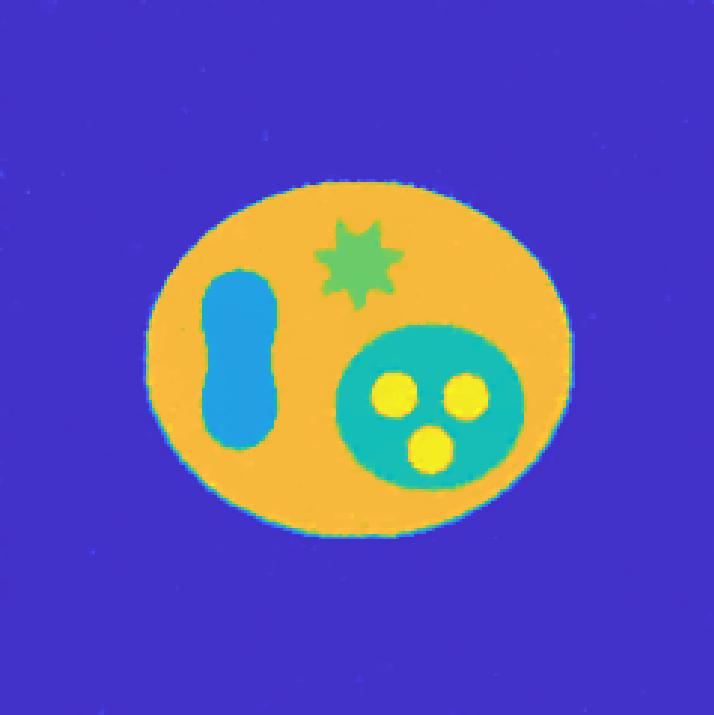};
  \end{axis}
\end{tikzpicture}
\caption{Reconstructions with moving axis (5\,\% Gaussian noise).
  Left: Filtered backpropagation (PSNR 21.38, SSIM 0.078).
  Center: CG Reconstruction (PSNR 23.75, SSIM 0.186). 
  Right: Primal-dual with TV and $\lambda=0.05$ (PSNR 34.14, SSIM 0.900). 
  \label{fig:cell3d-rec-noise}}
\end{figure}

\section{Phase retrieval}
\label{sec:pr}

In many practical applications of diffraction tomography,
only the intensities
\begin{equation*}
  \abs{\utot_t(x_1,x_2,\rM) },
  \quad (x_1, x_2) \in \R^2,\ t \in [0,T],
\end{equation*}
are available
since the measurement of the phase is technically very challenging and expensive.
The missing phase information of the complex-valued function $\utot_t$
thus has to be recovered in advance.
Especially because the inversion of the given magnitudes
equipped with an arbitrary phase
would yield a solution of $f$,
phase retrieval in diffraction tomography is highly ill posed.
Such behaviour is also well known for classical phase retrieval problems, 
where the ambiguities can be characterized analytically \cite{BeiPlo15,BenBeiEld17,BeiPlo20}.
Therefore,
we require a priori knowledge to obtain a meaningful reconstruction.
For our numerical experiments,
we exploit
that $f$ is nonnegative with bounded support,
and that $f$ often has a cartoon-like appearance,
e.g.,
for images of biological cells.
Unfortunately,
our phase retrieval problem in diffraction tomography fits
neither into the classical setting for instance studied by Gerchberg, Saxton, and Fienup \cite{GerSax72,Fie82}
nor into the modern approaches of PhaseLift \cite{CanEldStrVor13}
and its tensor-free variant \cite{BeiBre21},
since the connection between $\utot_t$ and $f$ is not linear but affine.

For phase retrieval in diffraction tomography,
we propose an all-at-one approach utilizing the \emph{hybrid input-output} (HIO) scheme,
see \cite[§\,6]{BeiQue22},
which is based on Fienup's HIO method \cite{Fie82}.
In a nutshell,
the main idea behind the iterative HIO approach
consists of computing the forward model,
projecting to the known magnitudes in the measurement domain,
solving the known-phase reconstruction problem with one of the above-mentioned methods
and incorporating the additional a priori information on $f$.
Instead of strictly enforcing the a priori constraints,
the HIO algorithm tries to correct the last input
to push the next output in direction of a feasible object.
In more details,
we apply the HIO method in Algorithm~\ref{alg:err-red},
where the operator $\bm D$ models the measurement process that maps the discretizations of $f$ to $\utot_t(\cdot,\cdot,\rM)$, $t\in[0,T]$, the sign function
\begin{equation*}
  \sgn(z) \coloneqq
  \begin{cases}
    z / |z|, & z \in\C\setminus \{0\}, \\
    1, & z = 0,
  \end{cases}
\end{equation*}
is applied element by element, 
and $\odot$ denotes the element-wise multiplication.

\begin{algorithm}[hbt]
  \KwIn{$\bm d = \abs{\bm D(\bm f)}$, $\beta \in (0,1]$, $\rs > 0$.}
  
  Initialize $\bm g^{(0)} \coloneqq \bm d$\;
  
  \For{$j=0,1,2,\dots,J_{\mathrm{HIO}}$}{
    $\bm f^{(j)} \coloneqq \bm D^{-1} \bm g^{(j)}$\;
    \For{$\bl\in\mathcal I$}{
      \eIf{$\| \bx_\bl \|_2 \le \rs$}
      {$\tilde f_\bl^{(j)}  \coloneqq \max\{f^{(j)}_{\bl} , 0 \}$\;}
      {$\tilde f_\bl^{(j)}  \coloneqq 0 $\;}
      \eIf{$f^{(j)}_{\bl} = \tilde f^{(j)}_{\bl}$}
      {$f^{(j+1/2)}_{\bl} \coloneqq f^{(j)}_{\bl}$\;}
      {$f^{(j+1/2)}_{\bl} \coloneqq f^{(j-1/2)}_{\bl} - \beta (f^{(j)}_{\bl} - \tilde f^{(j)}_{\bl})$\;}
    }
    $\bm g^{(j+1/2)} \coloneqq \bm D \bm f^{(j+1/2)}$\;
    $\bm g^{(j+1)} \coloneqq \bm d \odot \sgn (\bm g^{(j+1/2)})$\;
  }
  \KwOut{Approximate scattering potential $\bm f^{(J_{\mathrm{HIO}})}$.}
  \caption{Hybrid Input-Output Algorithm} 
  \label{alg:err-red}
\end{algorithm}

The crucial point of the proposed method is
that when computing $\bm D^{-1}$ we have to solve the known-phase reconstruction in each step,
which is performed via an iterative algorithm by itself.
Since both iterative methods show a slow convergence,
we made several numerical improvements to speed up computation.
For instance,
we restart the inner loop of the primal-dual algorithm with the parameters and dual variables
obtained in the previous outer step
and execute only a small number of inner iterations.
In this manner,
both iterative algorithms may be incorporated in each other.
Furthermore,
an initial reconstruction using CG as inner algorithm serves as starting point
for the more accurate but slower HIO with primal-dual iteration.
It is well known
that the performance of HIO may be improved
by utilizing a total variation denoising
between the iterates \cite{GauMohKha15,GauKha19}.
If the total variation regularization is used
during the backward transform,
this TV denoising is automatically incorporated into the HIO method.

\begin{figure}
  \centering
  \pgfmathsetmacro{\cmin} {-0.05} \pgfmathsetmacro{\cmax} {1.05}
  \begin{tikzpicture}
    \begin{axis}[
      width=\modelwidth,
      height=\modelwidth,
      enlargelimits=false,
      scale only axis,
      axis on top,
      colorbar,colorbar style={
        width=.15cm, xshift=-0.9em, 
        point meta min=\cmin,point meta max=\cmax,
      },
      ]
      \addplot graphics [
      xmin=\xmin, xmax=\xmax,  ymin=\xmin, ymax=\xmax,
      ] {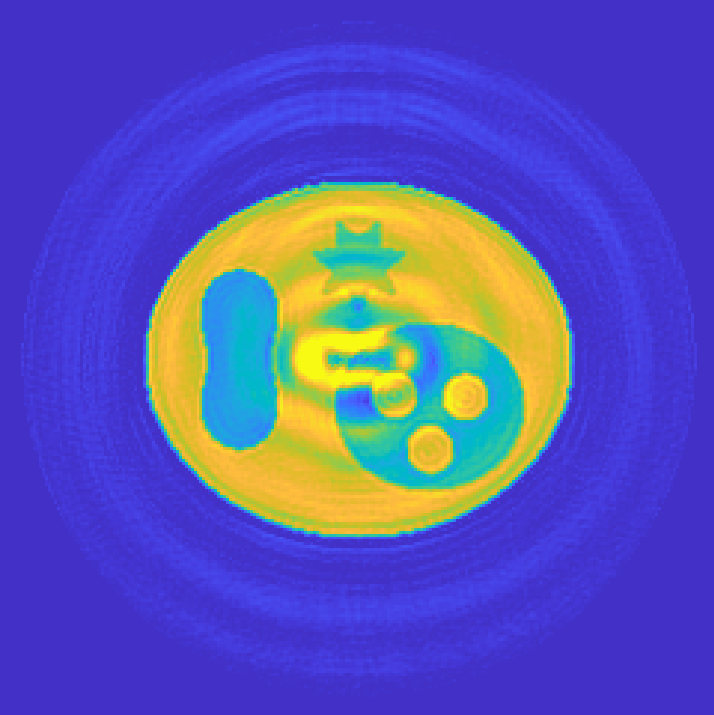};
    \end{axis}
  \end{tikzpicture}%
  \qquad
  \begin{tikzpicture}
    \begin{axis}[
      width=\modelwidth,
      height=\modelwidth,
      enlargelimits=false,
      scale only axis,
      axis on top,
      colorbar,colorbar style={
        width=.15cm, xshift=-0.9em, 
        point meta min=\cmin,point meta max=\cmax,
      },
      ]
      \addplot graphics [
      xmin=\xmin, xmax=\xmax,  ymin=\xmin, ymax=\xmax,
      ] {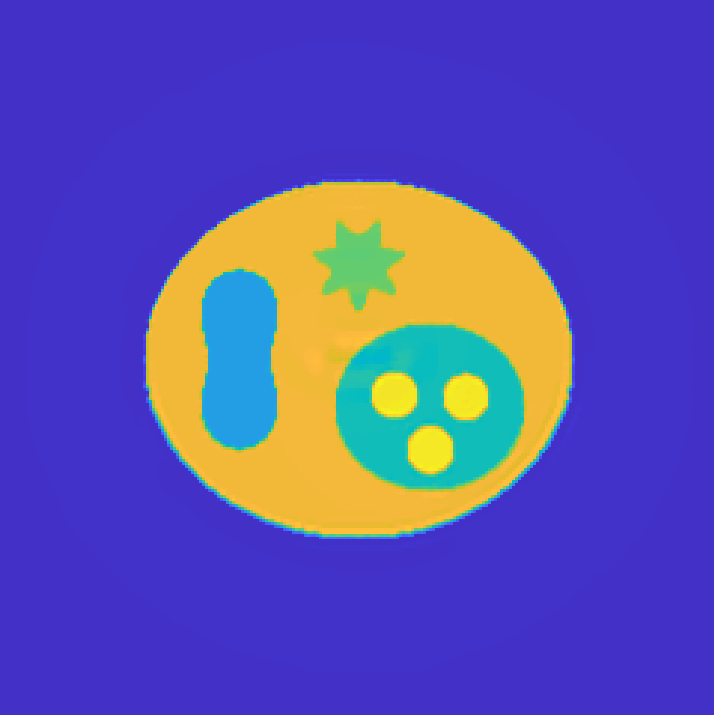};
    \end{axis}
  \end{tikzpicture}
  \caption{
    Phase retrieval (exact data).
    Left: HIO and CG method (PSNR 28.41, SSIM 0.697).
    Right: HIO primal-dual method (PSNR 34.58, SSIM 0.973). 
    \label{fig:cell3d-pr}}
\end{figure}

\begin{figure}
  \centering
  \pgfmathsetmacro{\cmin} {-0.05} \pgfmathsetmacro{\cmax} {1.05}
  \begin{tikzpicture}
    \begin{axis}[
      width=\modelwidth,
      height=\modelwidth,
      enlargelimits=false,
      scale only axis,
      axis on top,
      colorbar,colorbar style={
        width=.15cm, xshift=-0.9em, 
        point meta min=\cmin,point meta max=\cmax,
      },
      ]
      \addplot graphics [
      xmin=\xmin, xmax=\xmax,  ymin=\xmin, ymax=\xmax,
      ] {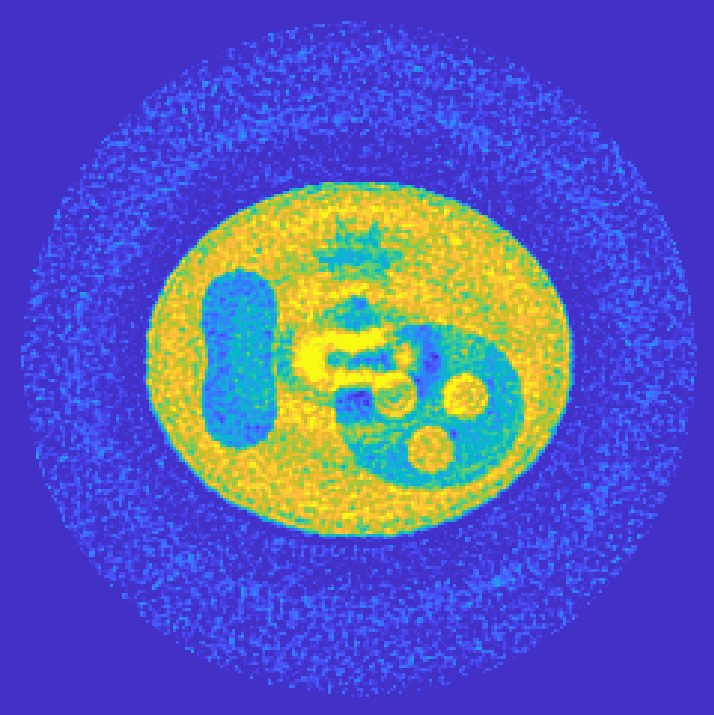};
    \end{axis}
  \end{tikzpicture}%
  \qquad
  \begin{tikzpicture}
    \begin{axis}[
      width=\modelwidth,
      height=\modelwidth,
      enlargelimits=false,
      scale only axis,
      axis on top,
      colorbar,colorbar style={
        width=.15cm, xshift=-0.9em, 
        point meta min=\cmin,point meta max=\cmax,
      },
      ]
      \addplot graphics [
      xmin=\xmin, xmax=\xmax,  ymin=\xmin, ymax=\xmax,
      ] {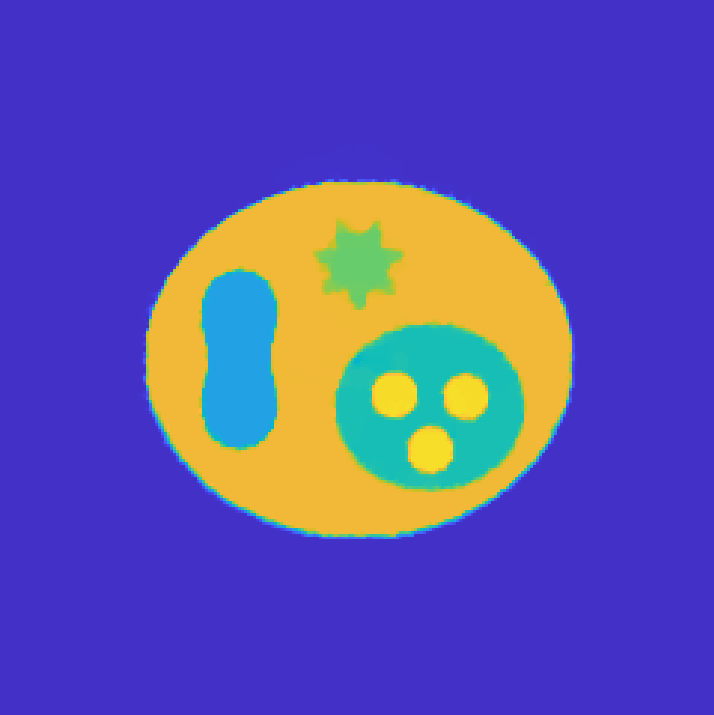};
    \end{axis}
  \end{tikzpicture}
  \caption{
    Phase retrieval (5\,\% Gaussian noise).
    Left: HIO and CG method (PSNR 23.47, SSIM 0.579).
    Right: HIO primal-dual method with $\lambda=0.05$ (PSNR 34.00, SSIM 0.991). 
    \label{fig:cell3d-pr-noise}}
\end{figure}

Numerical simulations indicate that phase retrieval in our diffraction tomography setup is indeed possible
utilizing the proposed HIO algorithm with the primal-dual TV regularization,
see Fig.~\ref{fig:cell3d-pr} and \ref{fig:cell3d-pr-noise}.
Here we use the same rotation axes \eqref{eq:axis} as before.
For the HIO with CG, we performed $J_{\mathrm{HIO}}=10$ outer and 5 inner (CG) iterations, where more iterations yield a noisier result.
For the HIO with primal-dual, we performed $J_{\mathrm{HIO}}=200$ outer and 5 inner iterations.
We note that, in contrast to the reconstructions of §\,\ref{sec:rec},
we incorporate additional a priori information via the support constraint $\rs=40$.
A Matlab implementation of the algorithms is available at \url{https://github.com/michaelquellmalz/FourierODT}.

\section{Conclusion}

We have compared reconstruction algorithms in diffraction tomography with a rigid object that undergoes a time-dependent rotation according to a moving rotation axis.
All three algorithms work well for exact data with an advantage for the iterative methods (CG and primal-dual).
In the presence of noise in the measured data, the TV-regularized primal dual shows superior results.
If only the intensity data $|\utot|$ is known, the all-at-once HIO algorithm combined with the TV regularization is still able to recover the object accurately.

\paragraph{Acknowledgement}
  R. Beinert gratefully acknowledges funding by the BMBF under the project “VI-Screen”
  (13N15754). M. Quellmalz gratefully acknowledges funding by the DFG (STE 571/19-1,
  project number 495365311) within SFB F68 (“Tomography Across the Scales”).

\vspace{\baselineskip}
%% The style of the following references should be used in all documents.

%\bibliography{abrv,references}
% \bibliographystyle{pamm}
% \bibliography{references}

\providecommand{\WileyBibTextsc}{}
\let\textsc\WileyBibTextsc
\providecommand{\othercit}{}
\providecommand{\jr}[1]{#1}
\providecommand{\etal}{~et~al.}

\end{document}